\documentstyle[12pt,amscd,amssymb]{amsart}

\newcommand{\exseq}[3]{0 \ar #1 \ar #2 \ar #3 \ar 0}
\newcommand{\seq}[3]{{#1}_{#2}, \ldots, {#1}_{#3}}

\newcommand{\M}{M_{\Sigma}}
\newcommand{\Y}{\Sigma \times {\Bbb S}^1}
\newcommand{\X}{\Sigma \times {\Bbb P}^1}
\newcommand{\Spz}{\text{Sp}\, (2g,{\Bbb Z})}

\newcommand{\surj}{\twoheadrightarrow}
\newcommand{\inc}{\hookrightarrow}
\newcommand{\ar}{\rightarrow}

\newcommand{\x}{\times}
\newcommand{\ox}{\otimes}
\newcommand{\iso}{\cong}
\newcommand{\isom}{\stackrel{\simeq}{\ar}}

\newcommand{\point}{\text{pt}}

\newcommand{\End}{\text{End}}
\newcommand{\Diff}{\text{Diff}}

\newcommand{\Gr}{\text{Gr}\,}

\newcommand{\coker}{\text{coker}}

\newcommand{\ch}{\text{ch}\:}
\newcommand{\Ext}{\text{Ext}}

\newcommand{\id}{\text{id}}

\newcommand{\cE}{{\cal E}}
\newcommand{\cF}{{\cal F}}

\newcommand{\cM}{{\cal M}}

\newcommand{\cN}{{\cal N}}
\newcommand{\cL}{{\cal L}}
\newcommand{\cO}{{\cal O}}

\newcommand{\cU}{{\cal U}}
\newcommand{\cV}{{\cal V}}

\renewcommand{\AA}{{\Bbb A}}
\newcommand{\CC}{{\Bbb C}}

\newcommand{\EE}{{\Bbb E}}

\newcommand{\PP}{{\Bbb P}}
\newcommand{\QQ}{{\Bbb Q}}

\newcommand{\SS}{{\Bbb S}}
\newcommand{\ZZ}{{\Bbb Z}}

\renewcommand{\a}{\alpha}
\renewcommand{\b}{\beta}

\newcommand{\g}{\gamma}

\newcommand{\f}{\epsilon}

\renewcommand{\l}{\lambda}

\newcommand{\s}{\sigma}

\renewcommand{\o}{\omega}
\newcommand{\p}{\phi}
\newcommand{\q}{\psi}

\newcommand{\z}{\zeta}

\renewcommand{\S}{\Sigma}

\renewcommand{\L}{\Lambda}
\renewcommand{\P}{\Phi}
\newcommand{\Q}{\Psi}

\newcommand{\frg}{{\frak g}}
\newcommand{\frM}{{\frak M}}

\theoremstyle{plain}
\begingroup
\theorembodyfont{\sl}
\newtheorem{thm}{Theorem}
\newtheorem{cor}[thm]{Corollary}
\newtheorem{lem}[thm]{Lemma}
\newtheorem{prop}[thm]{Proposition}

\endgroup

\theoremstyle{definition}

\theoremstyle{remark}
\newtheorem{rem}[thm]{Remark}
\newtheorem{ex}[thm]{Example}

\title[Quantum cohomology of moduli space of stable bundles]{Quantum cohomology
of
  the moduli space of stable bundles over a Riemann surface}

\author{Vicente Mu\~noz}
\address{Departamento de \'Algegra, Geometr\'{\i}a y Topolog\'{\i}a \\ Facultad
de Ciencias \\ Universidad de M\'alaga \\ 29071 M\'alaga \\ Spain}
\email{vmunoz@@agt.cie.uma.es}
\thanks{\hbox{$^*$}Supported by a grant from Ministerio de Educaci\'on y
Cultura \\ 1991 Mathematics Subject Classification. Primary: 57R57. Secondary:
14D20, 58D27, 53C23.}
\date{November, 1997}

\begin{document}

\maketitle

\begin{abstract}
  We determine the quantum cohomology of the moduli space $\M$ 
  of odd degree rank two
  stable vector bundles over a Riemann surface $\S$ of genus $g \geq 1$. 
  This work together with~\cite{Floer} complete the proof of the existence
  of an isomorphism $QH^*(\M) \iso HF^*(\Y)$.
\end{abstract}

\section{Introduction}
\label{sec:intro}

Let $\S$ be a Riemann surface of genus $g \geq 2$ and let $\M$ denote the
moduli space
of flat $SO(3)$-connections with nontrivial second Stiefel-Whitney class $w_2$.

This is a smooth symplectic manifold
of dimension $6g-6$. Alternatively, we can consider $\S$ as a smooth complex
curve of
genus $g$. Fix a line bundle $\L$ on $\S$ of degree $1$, then 
$\M$ is the moduli space of  rank two 
stable vector bundles on $\S$ with determinant $\L$, which is a smooth
complex variety of complex dimension $3g-3$. The symplectic deformation class 
of $\M$ 
only depends on $g$ and not on the particular complex structure on $\S$.

The manifold $X=\M$ is a positive symplectic manifold with $\pi_2(X)=\ZZ$.
For such a manifold $X$, its quantum cohomology,
$QH^*(X)$, is well-defined 
(see~\cite{Ruan}~\cite{RT}~\cite{McDuff}~\cite{Piunikhin}).
As vector spaces, $QH^*(X)=H^*(X)$ (rational coefficients are understood),
but the multiplicative structure is different. 
Let $A$ denote the positive generator of $\pi_2(X)$, i.e. the
generator such that the symplectic form evaluated on $A$ is positive.
Let $N=c_1(X)[A] \in \ZZ_{>0}$. Then there is a natural $\ZZ/2N\ZZ$-grading 
for $QH^*(X)$, which 
comes from reducing the $\ZZ$-grading of $H^*(X)$. 
(For the case $X=\M$, $N=2$, so $QH^*(\M)$ is $\ZZ/4\ZZ$-graded).
The ring structure of 
$QH^*(X)$, called quantum multiplication, 
is a deformation of the usual cup product for $H^*(X)$. 
For $\a \in H^p(X)$, $\b \in H^q(X)$, we define the quantum 
product of $\a$ and $\b$ as 
$$
  \a \cdot \b =\sum_{d \geq 0} \P_{dA}(\a, \b),
$$
where $\P_{dA}(\a, \b) \in H^{p+q-2Nd} (X)$ is given by 
$<\P_{dA}(\a, \b),\g>=\Q^X_{dA}(\a, \b, \g)$, the Gromov-Witten invariant, for
all
$\g \in H^{\dim X-p-q+2Nd} (X)$. One has $\P_0(\a, \b)=\a \cup \b$. 
The other terms are the quantum
correction terms and they all live in lower degree parts of the
cohomology groups. It is a fact~\cite{RT} that the quantum product gives an
associative and graded commutative ring structure. 

To define the Gromov-Witten invariant, 
let $J$ be a generic almost complex
structure compatible with the symplectic form. 
Then for every $2$-homology class $dA$, $d \in \ZZ$, 
there is a moduli space $\cM_{dA}$ of pseudoholomorphic rational 
curves (with respect to $J$) 
$f : \PP^1 \ar X$ with $f_*[\PP^1] = dA$.
Note that $\cM_0=X$ and that $\cM_{dA}$ is empty for $d<0$.
For $d \geq 0$,
the dimension of $\cM_{dA}$ is $\dim X +2Nd$. 
This moduli space
$\cM_{dA}$ admits a natural compactification, $\overline{\cM}_{dA}$, called
the
Gromov-Uhlenbeck compactification~\cite{Ruan}~\cite[section 3]{RT}.
Consider now $r \geq 3$ different points $\seq{P}{1}{r} \in {\PP}^1$. Then we
have defined an evaluation map $ev: \cM_{dA} \ar X^r$
by $f \mapsto (f(P_1), \ldots, f(P_r))$. This map extends to 
$\overline{\cM}_{dA}$ and its image, $ev(\overline{\cM}_{dA})$, is a
pseudo-cycle~\cite{RT}. So for $\a_i \in
H^{p_i}(\M)$, $1 \leq i \leq r$, with $p_1+ \cdots +p_r = \dim X + 2Nd$, we 
choose generic cycles $A_i$, $1 \leq i \leq r$,
representatives of their Poincar\'e duals, and set
\begin{equation}
\label{eqn:qu1}
  \Q_{dA}^X(\seq{\a}{1}{r}) = <A_1\x\cdots\x A_r, [ev(\overline{\cM}_{dA})]>= 
  \# ev_{P_1}^*(A_1) \cap \cdots \cap ev_{P_r}^*(A_r),
\end{equation}
where $\#$ denotes count of points (with signs) and
$ev_{P_i}: \cM_{dA} \ar X$, $f \mapsto f(P_i)$.
This is a well-defined number and independent of the particular cycles.
Also, as the manifold $X$ is positive, 
$\coker L_f=H^1(\PP^1, f^* c_1(X))=0$, for all $f \in \cM_{dA}$
(see~\cite{Ruan} for definition of $L_f$). By~\cite{Ruan}
the complex structure of $X$ is generic and we can use it to compute the
Gromov-Witten
invariants.

Also for $r \geq 2$, let $\a_i \in H^{p_i}(\M)$, $1 \leq i \leq r$, then
$$
  \a_1 \cdots \a_r =\sum_{d \geq 0} \P_{dA}(\a_1, \ldots ,\a_r),
$$ 
where the correction terms $\P_{dA}(\a_1, \ldots ,\a_r) \in
H^{p_1+\cdots +p_r - 2 N d}(X)$ are determined by 
$<\P_{dA}(\a_1, \ldots ,\a_r), \g> = \Q_{dA}^X(\a_1, \ldots ,\a_r, \g)$, for
any 
$\g \in H^{\dim X +2Nd-(p_1+\cdots +p_r)}(X)$.

Returning to our manifold $X=\M$, there is a classical conjecture relating
the quantum cohomology $QH^*(\M)$ and the instanton Floer cohomology of
the three manifold
$\Y$, $HF^*(\Y)$ (see~\cite{Floer}).
In~\cite{Vafa} a presentation of $QH^*(\M)$ was given using physical 
methods, and in~\cite{Floer} it was proved that such a presentation 
was a presentation of $HF^*(\Y)$ indeed. 
Here we determine a presentation of $QH^*(\M)$ and prove the isomorphism
$QH^*(\M) \iso HF^*(\Y)$.

Siebert and Tian have an alternative 
program~\cite{Siebert} to find the presentation
of $QH^*(\M)$, which goes through proving a recursion formula for the
Gromov-Witten invariants of $\M$ in terms of the genus $g$.

The paper is organised as follows. 
In section~\ref{sec:ordinary} we review the 
ordinary cohomology ring of $\M$. In section~\ref{sec:curves} the moduli space
of
lines (rational curves representing $A$) in $\M$ is described. This makes
possible to compute the Gromov-Witten 
invariants $\Q_A^{\M}$, which determine
the first quantum correction terms of the quantum products in $QH^*(\M)$.
Section~\ref{sec:GW-inv} is devoted to this task. 
In~\cite{D1} Donaldson uses this information 
alone to determine $QH^*(\M)$ in the case of
genus $g=2$. It is somehow natural to expect that this idea can be developed in
the
general case $g \geq 3$.
In section~\ref{sec:quantum} we give an explicit presentation of
$QH^*(\M)$ for $g \geq 3$
(theorem~\ref{thm:main}), concluding the proof of $QH^*(\M) \iso HF^*(\Y)$
(corollary~\ref{cor:21}).
The two main ingredients that we make use of are the $\Spz$-decomposition of
$H^*(\M)$
under the action of the mapping class group (not ignoring the non-invariant
part
as it was customary) and a recursion similar to that in~\cite{Siebert} 
(lemma~\ref{lem:17}). The difference with~\cite{Siebert} lies in the fact that
we 
fix the genus, so that we do not need to compare the Gromov-Witten invariants
for 
moduli spaces of Riemann surfaces of different genus.
Finally in section~\ref{sec:6} we discuss the cases $g=1$ and $g=2$, which are
slightly
different.

\noindent {\em Acknowledgements:\/} The author is greatly indebted to 
the two sources of inspirations who have made possible this work:
Prof.\ Simon Donaldson, for his encouragement and invaluable help,
and Bernd Siebert and Gang Tian, 
for helpful conversations and for providing with a copy of their 
preprint~\cite{Siebert}, which was very enlightening.

\section{Classical cohomology ring of $\M$}
\label{sec:ordinary}

Let us recall the known description of the homology of
$\M$~\cite{King}~\cite{ST}~\cite{Floer}.
Let $\cU \ar \S \x \M$ be the universal bundle and consider the K\"unneth
decomposition of
\begin{equation}
   c_2(\End_0 \, \cU)=2 [\S] \ox \a + 4 \q -\b,
\label{eqn:qu2}
\end{equation}
with $\q=\sum \g_i \ox \q_i$, where 
$\{\seq{\g}{1}{2g}\}$ is a symplectic basis of 
$H^1(\S;\ZZ)$ with $\g_i \g_{i+g}=[\S]$ for $1 \leq i \leq g$ (also
$\{\g^{\#}_i\}$ will denote the dual basis for $H_1(\S;\ZZ)$). 
Here we can suppose without loss of generality that $c_1(\cU)=\L + \a$
(see~\cite{ST}).
In terms of the map $\mu: H_*(\S) \to H^{4-*}(\M)$,
given by $\mu(a)= -{1 \over 4} \,
p_1(\frg_{\cU}) /a $ (here $\frg_{\cU} \to \S \x \M$ is the associated
universal $SO(3)$-bundle, and $p_1(\frg_{\cU}) \in H^4(\S\x\M)$ its first
Pontrjagin class), we have
$$
   \left\{ \begin{array}{l} \a= 2\, \mu(\S) \in H^2
   \\ \q_i= \mu (\g_i^{\#}) \in H^3, \qquad 1\leq i \leq 2g
   \\ \b= - 4 \, \mu(x) \in H^4       
    \end{array} \right.
$$
where $x \in H_0(\S)$ is the class of the point, and $H^i=H^i(\M)$.
These elements generate $H^*(\M)$ as a 
ring~\cite{King}~\cite{Thaddeus}, and $\a$ is the positive generator of
$H^2(\M;\ZZ)$. We can rephrase this as saying that there exists an 
epimorphism 
\begin{equation}
  \AA(\S)= \QQ[\a,\b ]\ox \L(\seq{\q}{1}{2g}) \surj H^*(\M)
\label{eqn:qu3}
\end{equation}
(the notation $\AA(\S)$ follows that of Kronheimer and Mrowka~\cite{KM}, 
although it is slightly different).
Recall that $\deg(\a)=2$, $\deg(\b)=4$ and $\deg(\q_i)=3$.

The mapping class group $\Diff(\S)$ acts on $H^*(\M)$, with the action
factoring
through the action of $\Spz$ on $\{\q_i\}$. 
The invariant part, $H_I^*(\M)$, is generated by
$\a$, $\b$ and $\g=-2 \sum_{i=0}^g \q_i\q_{i+g}$. Then
there is an epimorphism
\begin{equation}
  \QQ[\a,\b,\g ]  \surj H^*_I(\M)
\label{eqn:qu4}
\end{equation}
which allows us to write
$$
   H_I^*(\M)= \QQ [\a, \b, \g]/I_g,
$$
where $I_g$ is the ideal of relations satisfied by $\a$, $\b$ and $\g$. 
{}From~\cite{ST}, a
basis for $H_I^*(\M)$ is given by the monomials $\a^a\b^b\g^c$, with
$a+b+c<g$.
For
$0 \leq k \leq g$, the primitive component of $\L^k H^3$ is 
$$
   \L_0^k H^3 = \ker (\g^{g-k+1} : \L^k H^3 \ar \L^{2g-k+2} H^3).
$$
The spaces $\L^k_0 H^3$ are irreducible $\Spz$-modules, i.e. the transforms
of any nonzero element of $\L^k_0 H^3$ under $\Spz$ generate the whole of it. 
The description of the ideals $I_g$ and the cohomology ring
$H^*(\M)$ is given in the following

\begin{prop}[\cite{ST}~\cite{King}]
\label{prop:1}
  Define $q^1_0=1$, $q^2_0=0$, $q^3_0=0$ and then recursively, for all 
  $r \geq 1$,
  $$
   \left\{ \begin{array}{l} q_{r+1}^1 = \a q_r^1 + r^2 q_r^2
   \\ q_{r+1}^2 = (\b+(-1)^{r+1}8) q_r^1 + {2r \over r+1}  q_r^3
   \\ q_{r+1}^3 =  \g  q_r^1
    \end{array} \right.
  $$
  Then $I_g=(q^1_g,q^2_g,q^3_g) \subset \QQ[\a,\b,\g]$, for all $g \geq 1$. 
  Note that $\deg(q_g^1)=2g$, $\deg(q_r^2)=2g+2$ and $\deg(q_g^3)=2g+4$.
  Moreover the $\Spz$-decomposition of $H^*(\M)$ is 
\begin{equation}
   H^*(\M)= \bigoplus_{k=0}^{g-1} \L_0^k H^3 \ox \QQ [\a, \b, \g]/I_{g-k}.
\label{eqn:antes}
\end{equation}
\end{prop}

This proposition allows us to find a basis for $H^*(\M)$ as follows. 
Let $\{x_i^{(k)}\}_{i \in B_k}$ be a 
basis of $\L_0^k H^3$, $0 \leq k \leq g-1$. Then 
\begin{equation} 
  \{x^{(k)}_i\a^a\b^b\g^c/ k=0,1,\ldots, g-1, \> a+b+c < g-k, \> i \in B_k \}
\label{eqn:antes2}
\end{equation}
is a basis for $H^*(\M)$. If we set
\begin{equation} 
  x^{(k)}_0= \q_1 \q_2 \cdots \q_k \in\L_0^k H^3,
\label{eqn:antes3}
\end{equation} 
then proposition~\ref{prop:1} says that a complete
set of relations satisfied in $H^*(\M)$
are $x_0^{(k)} q^i_{g-k}$, $i=1,2,3$, $0 \leq k
\leq g$, and the $\Spz$ transforms of these.

\section{Holomorphic lines in $\M$}
\label{sec:curves}

In order to compute the Gromov-Witten invariants $\Q_A^{\M}$, we need
to describe the space of lines, i.e. rational 
curves in $\M$ representing the generator $A \in H_2(\M;\ZZ)$,
$$
 \cM_A=\{f: \PP^1 \ar \M / \text{$f$ holomorphic, } f_*[\PP^1]=A\}.
$$

Let us fix some notation. 
Let $J$ denote the Jacobian variety of $\S$ parametrising line bundles of
degree $0$ and
let $\cL \ar \S\x J$ be the universal line bundle. If $\{\g_i\}$ is
the basis of $H^1(\S)$ introduced in section~\ref{sec:ordinary}
then $c_1({\cL}) = \sum \g_i \ox \p_i \in H^1(\S) \ox H^1(J)$, where 
$\{\p_i\}$ is a symplectic basis for $H^1(J)$. Thus $c_1({\cL})^2 =
-2 [\S] \ox \o \in H^2(\S) \ox H^2(J)$, where $\o=\sum_{i=1}^g \p_i \wedge
\p_{i+g}$ is the natural symplectic form for $J$.

Consider now the algebraic surface $S=\X$. It has irregularity $q=g \geq 2$, 
geometric genus $p_g=0$ and canonical bundle $K \equiv -2\S+(2g-2){\PP}^1$.
Recall that $\L$ is a fixed line bundle of degree $1$ on $\S$.
Fix the line bundle $L=\L \ox \cO_{\PP^1}(1)$ on $S$ (we omit all
 pull-backs)
with $c_1= c_1(L) \equiv \PP^1+\S$, 
and put $c_2=1$. 
The ample cone of $S$ is $\{a\PP^1 +b \S \, / \, a,b >0 \}$.
Let $H_0$ be a polarisation close to $\PP^1$ in the ample cone
and $H$ be a polarisation close to $\S$, i.e.
$H= \S +t\PP^1$ with $t$ small.
We wish to study the moduli space $\frM= \frM_H(c_1,c_2)$ of $H$-stable
bundles
over $S$ with Chern classes $c_1$ and $c_2$. 

\begin{prop}
\label{prop:2}
  $\frM$ can be described as a bundle  $\PP^{2g-1} \ar \frM=
\PP(\cE_{\z}^{\vee})
  \ar J$, where $\cE_{\z}$ is a bundle on $J$ with $\ch \cE_{\z} = 2g + 8 \o$.
  So $\frM$ is compact, smooth and of the expected dimension $6g-2$.
  The universal bundle $\cV \ar S \x \frM$ is given by
$$
  0 \ar \cO_{\PP^1}(1) \ox \cL \ox \l \ar \cV \ar \L \ox \cL^{-1} \ar 0,
$$
  where $\l$ is the tautological line bundle for $\frM$.
\end{prop}

\begin{pf}
  For the polarisation $H_0$, the moduli space of $H_0$-stable 
  bundles with Chern classes
  $c_1, c_2$ is empty by~\cite{Qin}. Now for $p_1=-4c_2+c_1^2=-2$ 
  there is only one wall,
  determined by $\z \equiv -\PP^1 +\S$
  (here we fix $\z= 2\S- L=\S-c_1(\L)$ as a divisor), 
  so the moduli space of $H$-stable bundles with 
  Chern classes $c_1, c_2$ is obtained by crossing the wall as described
  in~\cite{wall}.
  First, note that the results in~\cite{wall} use the hypothesis of $-K$
  being 
  effective, but the arguments work equally well with the weaker assumption
  of $\z$ being a good wall~\cite[remark 1]{wall} (see also~\cite{Gottsche}
  for
  the case of $q=0$). In our case, $\z \equiv -\PP^1+\S$ is a good wall
  (i.e. 
  $\pm \z+K$ are both not effective) with $l_{\z}=0$. Now with the notations
  of~\cite{wall}, $F$ is a divisor such that $2F-L \equiv \z$, e.g. $F=\S$.
  Also
  $\cF \ar S \x J$ is the universal bundle parametrising divisors
  homologically
  equivalent to $F$, i.e. $\cF=\cL \ox \cO_{\PP^1}(1)$. Let
  $\pi: S \x J \ar J$ be
  the projection. Then $\frM=E_{\z}= \PP(\cE_{\z}^{\vee})$, where
$$
  \cE_{\z} = {\cE}\text{xt}^1_{\pi}(\cO(L-\cF),\cO(\cF))=  R^1\pi_* 
  (\cO(\z )\ox \cL^2).
$$
  Actually $\frM$ is exactly the set
  of bundles $E$ that can be written as extensions
$$
  \exseq{\cO_{\PP^1}(1) \ox L}{E}{\L \ox L^{-1}}
$$
  for a line bundle $L$ of degree $0$.
  The Chern character is computed in~\cite[section 3]{wall} to be $\ch \cE_{\z}
= 2g +
  e_{K-2\z}$,
  where $e_{\a}= -2 (\PP^1 \cdot \a) \o$ (the class $\S$ defined
  in~\cite[lemma 11]{wall}
  is $\PP^1$ in our case).
  Finally, the description of the universal bundle follows
  from~\cite[theorem 10]{wall}.
\end{pf}

\begin{prop}
\label{prop:3}
  There is a well defined map $\cM_A \ar \frM$.
\end{prop}

\begin{pf}
  Every line $f:{\PP}^1 \ar \M$ gives a bundle $E=(\id_{\S}\x f)^*\cU$ over
  $\X$
  by pulling-back the universal bundle ${\cU} \ar \S \x \M$. Then for any 
  $t \in \PP^1$, the bundle $E|_{\S \x t}$ is defined by $f(t)$.
  Now, by
  equation~\eqref{eqn:qu2}, $p_1(E)= p_1(\cU)[\S \x A]= -2\a [A] =-2$.
  Since $c_1(E) = (\id_{\S}\x f)^*c_1(\cU)=\L +\S$, it must be $c_2=1$.
  To see that $E$ is $H$-stable, consider any sub-line bundle $L \inc E$
  with $c_1(L) \equiv a\PP^1 +b\S$. Restricting to any $\S \x t \subset \X$ 
  and using the
  stability of $E|_{\S \x t}$, one gets $a \leq 0$. Then $c_1(L) \cdot \S < 
  {c_1(E)\cdot \S \over 2}$, which yields the $H$-stability of $E$
  (recall that $H$ is close to $\S$). So $E \in \frM$.
\end{pf}

Now define $N$ as the set of extensions on $\S$ of the form 
\begin{equation} 
  0 \ar L \ar E \ar \L \ox L^{-1} \ar 0,
\label{eqn:qu5}
\end{equation}
for $L$ a line bundle of degree $0$. Then the groups 
$\Ext^1( \L \ox L^{-1}, L) =H^1(L^2 \ox
\L^{-1})= H^0(L^{-2} \ox \L \ox K)$ are of constant dimension $g$. Moreover 
$H^0(L^{2} \ox \L^{-1})=0$, so the moduli space $N$ which
parametrises extensions like~\eqref{eqn:qu5} is given as 
$N=\PP(\cE^{\vee})$,
where $\cE= {\cE}{\hbox{xt}}^1_p(\L \ox {\cL}^{-1},
{\cL})=R^1 p_*({\cL}^{2} \ox
\L^{-1})$, $p:\S\x J \ar J$ the projection.
Then we have a fibration ${\PP}^{g-1} \ar N=\PP(\cE^{\vee}) \ar J$. 
The Chern character of $\cE$ is
\begin{equation} 
\begin{array}{rcl}
\ch({\cal E}) &=& \ch(R^1 p_*({\cL}^{2}\ox \L ^{-1})) =
             - \ch(p_{!}( {\cL}^2 \ox \L ^{-1}) ) = \\
             &=& -p_*( (\ch\,{\cL})^2 \> (\ch\,\L)^{-1} \> \text{Todd} \> 
              T_{\S}) = \\
             &=& -p_*((1+c_1({\cL})+{1 \over
               2}c_1({\cL})^2)^2(1-\L)(1-{1 \over 2}K)) = \\
             &=& -p_*(1-{1 \over 2}K+2c_1({\cL}) -4 \o
               \ox [\S]-\L) = g+4\o.
\end{array} \label{eqn:qu6}
\end{equation}
It is easy to check that 
all the bundles in $N$ are stable, so there is a well-defined map 
$$ 
 i: N \ar \M.
$$

Now we wish to construct the space of lines in $N$.
Note that $\pi_2(N)=\pi_2(\PP^{g-1})=\ZZ$, 
as there are no rational curves in $J$. Let
$L \in \pi_2(N)$ be the positive generator. We want to 
describe 
$$
   \cN_L=\{f: \PP^1 \ar N / \text{$f$ holomorphic, } f_*[\PP^1]=L\}.
$$
For the projective space $\PP^n$, the space $H_1$ of lines in $\PP^n$ 
is the set of algebraic maps $f: {\PP}^1
\ar {\PP}^n$ of degree $1$. Such an $f$ has the form
$f[  x_0, x_1 ]=[x_0 u_0 +x_1 u_1]$, $[x_0,x_1] \in \PP^1$, 
where $u_0$, $u_1$ are linearly independent vectors 
in $\CC^{n+1}$. So
$$
  H_1= \PP(\{(u_0, u_1)/ \text{$u_0$, $u_1$ are linearly independent} \})
  \subset
  \PP((\CC^{\,n+1} \oplus \CC^{\,n+1})^{\vee}) =\PP^{2n+1}.
$$
The complement of $H_1$ is the image of $\PP^n \x \PP^1 \inc \PP^{2n+1}$,
$([u] , [x_0, x_1]) \mapsto [x_0 u, x_1 u]$, which is a smooth
$n$-codimensional
algebraic subvariety.
So $\cN_L$ can be described as the fibration
\begin{equation}
  \begin{array}{ccccc}
     H_1     &\ar &  \cN_L   & \ar& J \\
  \bigcap    &    & \bigcap  &    & \| \\  
  \PP^{2g-1} &\ar & \PP((\cE \oplus \cE)^{\vee}) & \ar & J
  \end{array}
\label{eqn:qu7}
\end{equation}

\begin{rem}
\label{rem:4}
Note that $\cE_{\z}= R^1\pi_* (\cO(\z )\ox \cL^2)= 
R^1\pi_*(\cO_{\PP^1}(1) \ox {\cL}^{2} \ox
\L^{-1})= H^0 (\cO_{\PP^1}(1)) \ox R^1p_*( {\cL}^{2} \ox \L^{-1})
= H^0 (\cO_{\PP^1}(1)) \ox \cE \iso \cE \oplus \cE$. So $\frM=\PP((\cE \oplus 
\cE)^{\vee})$, canonically.
\end{rem}

\begin{prop}
\label{prop:5}
The map $i:N \ar \M$ induces a map $i_*: \cN_L \ar \cM_A$. The composition
$\cN_L \ar \cM_A \ar \frM$ is the natural inclusion of~\eqref{eqn:qu7}.
\end{prop}

\begin{pf}
The first assertion is clear as $i$ is a holomorphic map.
For the second, consider the universal sheaf on $\S\x N$,
\begin{equation}
  0 \ar {\cL} \ox U \ar {\EE} \ar  \L \ox \cL^{-1} \ar 0,
\label{eqn:qu8}
\end{equation}
where $U=\cO_N(1)$ is the tautological bundle of
the fibre bundle $\PP^{g-1} \ar N \ar J$. 
Any element in $\cN_L$ is a line 
${\PP}^1 \inc N$, which must lie inside a single fibre $\PP^{g-1}$.
Restricting~\eqref{eqn:qu8} to this line, 
we have an extension 
$$
   0 \ar L \ox {\cO}_{{\PP}^1}(1) \ar E \ar \L \ox L^{-1}  \ar 0
$$
on $S=\X$, which is the image of the given element in $\frM$  
(here $L$ is the line bundle corresponding to the fibre 
in which ${\PP}^1$ sits). Now it is easy to check that the
map $\cN_L \ar \frM$ is the inclusion of~\eqref{eqn:qu7}.
\end{pf}

\begin{cor}
\label{cor:6}
$i_*$ is an isomorphism.
\end{cor}

\begin{pf}
By proposition~\ref{prop:5}, 
$i_*$ has to be an open immersion. The group $PGL(2,\CC)$ acts on both
spaces $\cN_L$ and $\cM_A$, and $i_*$ is equivariant.
The quotient $\cN_L/PGL(2,\CC)$ is compact, being
a fibration over the Jacobian with all the fibres the Grassmannian
$\Gr(\CC^{\,2},
\CC^{\,g-1})$, hence irreducible. As a consequence $i_*$ is an isomorphism.
\end{pf}

\begin{rem}
\label{rem:7}
Notice that the lines in $\M$ are all contained in the image of $N$,
which is of dimension $4g-2$ against $6g-6=\dim \M$.
They do not fill all of $\M$ as one would naively expect. 
\end{rem}

\section{Computation of $\Q_A^{\M}$}
\label{sec:GW-inv}

The manifold $N$ is positive with $\pi_2(N)=\ZZ$ and $L \in \pi_2(N)$ is
the positive generator. Under the map $i: N \ar \M$, 
we have $i_* L=A$. Now $\dim N= 4g-2$ and $c_1(N)[L]=c_1(\PP^{g-1})[L]=g$. So
quantum cohomology of $N$, $QH^*(N)$, is well-defined and
$\ZZ/2g\ZZ$-graded.
{}From corollary~\ref{cor:6}, it is straightforward to prove

\begin{lem}
\label{lem:8}
  For any $\a_i \in H^{p_i}(\M)$, $1 \leq i \leq r$, such that  $p_1 + \cdots +
p_r
  =6g-2$, it is
  $\Q_A^{\M}(\seq{\a}{1}{r})= \Q_L^N(\seq{i^*\a}{1}{r})$. $\quad \Box$
\end{lem}

It is therefore important to know the Gromov-Witten invariants of $N$, i.e.
its quantum cohomology.
{}From the universal bundle~\eqref{eqn:qu8}, we can read the first Pontrjagin
class $p_1( {\frg}_{\EE}) = -8  [\S]\ox\o + h^2- 2 [\S]\ox h
 + 4h\cdot c_1({\cL}) \in H^4(\S\x N)$, where $h=c_1(U)$ is the hyperplane
class. So on $N$ we have
\begin{equation}
  \left\{ \begin{array}{l}
  \a=2 \mu(\S) = 4\o + h \\
  \q_i=\mu(\g_i^{\#}) = - h \cdot \p_i \\
  \b= -4 \mu(x) =  h^2 
  \end{array} \right. 
\label{eqn:qu9}
\end{equation}
Let us remark here that $h^2$ denotes ordinary cup product in $H^*(N)$,
a fact which will prove useful later. Now
let us compute the quantum cohomology ring of $N$. The cohomology of $J$ is
$H^*(J)= \L H_1$, where $H_1=H_1(\S)$. Now the fibre bundle 
description $\PP^{2g-1} \ar N=\PP(\cE^{\vee}) \ar J$
implies that the usual cohomology of $N$ is
$H^*(N)= \L H_1 [h]/<h^g+c_1h^{g-1}+ \cdots +c_g=0>$, where  $c_i=c_i(\cE)={4^i
 \over i!} \o^i$, from~\eqref{eqn:qu6}.
As the quantum cohomology has the same generators as the
usual cohomology and the relations are a deformation
of the usual relations~\cite{ST2}, it must be $h^g+c_1h^{g-1}+ \cdots +c_g=r$
in $QH^*(N)$, with
$r \in \QQ$. As in~\cite[example 8.5]{RT}, $r$ can be computed to be $1$. So
\begin{equation}
  QH^*(N)= \L H_1 [h]/<h^g+c_1h^{g-1}+ \cdots +c_g=1> .
\label{eqn:qu10}
\end{equation}

\begin{lem}
\label{lem:9}
  For any $s \in H^{2g-2i}(J)$, $0 \leq i \leq g$, denote by $s \in
  H^{2g-2i}(N)$ its pull-back to $N$ under the natural projection. 
  Then the quantum product $h^{2g-1+i} s$
  in $QH^*(N)$ has component in $H^{4g-2}(N)$ equal to ${(-8)^i\over i!}
  \o^i \wedge s$ (the natural isomorphism $H^{4g-2}(N) \iso H^{2g}(J)$ is
  understood).
\end{lem}

\begin{pf}
  First note that for $s_1,s_2 \in H^*(J)$ such that their 
  cup product in $J$ is $s_1s_2=0$, then the quantum product $s_1s_2 \in
  QH^*(N)$ vanishes. This is so since every rational line in $N$ is 
  contained in a fibre of $\PP^{2g-1} \ar J \ar N$.

  Next recall that $h^{g-1+i}s$ has component in $H^{4g-2}(N)$ equal to
  $s_i(\cE) \wedge s={(-4)^i\over i!} \o^i \wedge s$.
  Then multiply the standard relation~\eqref{eqn:qu10}
  by $h^{g-1+i}s$ and work by induction on $i$. For $i=0$ we get
  $h^{2g-1}s=h^{g-1}s$ and the assertion is obvious. For $i>0$, 
$$
  h^{2g-1+i}s +h^{2g-2+i}c_1s+ \cdots +h^{2g-1}c_is=h^{g-1+i}s.
$$
  So the component of $h^{2g-1+i}s$ in $H^{4g-2}(N)$ is 
$$
  - \sum_{j=1}^i  {(-8)^{i-j}\over (i-j)!} \o^{i-j} c_j s +
  {(-4)^i\over i!} \o^i s =
  {(-8)^i\over i!} \o^i s
  - \sum_{j=0}^i  {(-8)^{i-j}\over (i-j)!}{4^j\over j!} \o^i s
  + {(-4)^i\over i!} \o^i s =   {(-8)^i\over i!} \o^i s.
$$
\end{pf}

\begin{lem}
\label{lem:10} 
Suppose $g>2$. Let $\a^a\b^b\q_{i_1}\cdots\q_{i_r} \in \AA(\S)$ 
have degree $6g-2$. Then 
$$
  \Q^N_L(\a, \stackrel{(a)}{\ldots}, \a, \b ,
  \stackrel{(b)}{\ldots} ,\b,\q_{i_1}, \ldots, \q_{i_r})=
  <(4\o + X)^a(X^2)^b \p_{i_1}\cdots\p_{i_r} X^r, [J]>,
$$
  evaluated on $J$, where $X^{2g-1+i}={(-8)^i  \over i!} \o^i \in H^*(J)$.
\end{lem}

\begin{pf}
  By definition the left hand side is the component in $H^{4g-2}(N)$
  of the quantum product $\a^a\b^b\q_{i_1}\cdots\q_{i_r} \in QH^*(N)$. 
  From~\eqref{eqn:qu9}, this quantum product is $(4\o + h)^a (h^2)^b
  (-h\p_{i_1})\cdots(-h\p_{i_r})$, 
  upon noting that when $g>2$, $\b=h^2$  
  as a quantum product as there are no quantum corrections 
  because of the degree.
  Note that $r$ is even, so the statement of the lemma follows from 
  lemma~\ref{lem:9}.
\end{pf}

Now we are in the position of relating the Gromov-Witten invariants 
$\Q_A^{\M}$ with the Donaldson invariants for $S=\S\x \PP^1$ 
(for definition of Donaldson invariants see~\cite{DK}~\cite{KM}).

\begin{thm}
\label{thm:11}
  Suppose $g>2$. Let $\a^a\b^b\q_{i_1}\cdots\q_{i_r} \in \AA(\S)$ have 
  degree $6g-2$. Then 
$$
  \Q^{\M}_A(\a, \stackrel{(a)}{\ldots}, \a, \b ,
  \stackrel{(b)}{\ldots} ,\b,\q_{i_1}, \ldots, \q_{i_r})=
  (-1)^{g-1} D^{c_1}_{S, H} ((2\S)^a(-4 \point)^b \g_{i_1}^{\#}
  \cdots\g_{i_r}^{\#}),
$$
  where $D^{c_1}_{S, H}$ stands for
  the Donaldson invariant of $S=\S\x\PP^1$ with 
  $w=c_1$ and polarisation $H$.
\end{thm}

\begin{pf}
  By definition, the right hand side is $\f_S(c_1)
  <\a^a\b^b\q_{i_1}\cdots\q_{i_r}, [\frM]>$, where 
  $\a =2\mu(\S) \in H^2(\frM)$, $\b =-4\mu(x) \in H^4(\frM)$,
  $\q_i =\mu(\g_i^{\#}) \in H^3(\frM)$.
  Here the factor $\f_S(c_1)=(-1)^{K_S c_1 +c_1^2
  \over 2}=(-1)^{g-1}$ compares the complex
  orientation of $\frM$ and its natural orientation as a moduli space of
  anti-self-dual connections~\cite{DK}. 
  By~\cite[theorem 10]{wall}, this is worked out to be 
  $(-1)^{g-1} <(4\o + X)^a(X^2)^b \p_{i_1}\cdots\p_{i_r} X^r,[J]>$, where 
  $X^{2g-1+i}=s_i(\cE_{\z})={(-8)^i  \over i!} \o^i$. Thus the 
  theorem follows from 
  lemmas~\ref{lem:8} and~\ref{lem:10}.
\end{pf}

\begin{rem}
\label{rem:12}
The formula in theorem~\ref{thm:11} is not right for $g=2$, as in
such case, the quantum product $h^2 \in QH^*(N)$ differs from $\b$ by
a quantum correction.
\end{rem}

\begin{rem}
\label{rem:13}
  Suppose $g \geq 2$ and let 
  $\a^a\b^b\q_{i_1}\cdots\q_{i_r} \in \AA(\S)$ have 
  degree $6g-6$. Then 
\begin{eqnarray*}
  \Q^{\M}_0(\a, \stackrel{(a)}{\ldots}, \a, \b ,
  \stackrel{(b)}{\ldots} ,\b,\q_{i_1}, \ldots, \q_{i_r}) &=&
  \f_S(\PP^1) < \a^a\b^b\q_{i_1}\cdots\q_{i_r},[\M]> = \\
  &=& - D^{\PP^1}_{S, H} ((2\S)^a(-4 \point)^b
  \g_{i_1}^{\#} \cdots\g_{i_r}^{\#}),
\end{eqnarray*}
as the moduli space of anti-self-dual connections on $S$ of dimension
$6g-6$ is $\M$.
\end{rem}

\section{Quantum cohomology of $\M$}
\label{sec:quantum}

It is natural to ask to what extent the first quantum correction 
determines the
full structure of the quantum cohomology of $\M$. In~\cite{D1}, Donaldson
finds the first quantum correction for $\M$ when the genus of $\S$ is 
$g=2$ and proves that this is enough to find the quantum product. 
Now it is our intention to show how the Gromov-Witten invariants $\Q^{\M}_A$
determine completely $QH^*(\M)$. First we check an interesting fact.

\begin{lem}
\label{lem:14}
  Let $g \geq 3$. Then $\g =-2 \sum \q_i\q_{i+g}$ 
  as elements in $QH^*(\M)$ (i.e.
  using the quantum product in the right hand side).
\end{lem}

\begin{pf}
  Let $\hat \g =-2 \sum \q_i\q_{i+g} \in QH^*(\M)$. 
  In principle, it is $\hat \g=\g + s\a$, for
  some $s \in \QQ$. Let us show that $s=0$. Multiplying by $\a^{3g-4}$, we
have
  $\hat \g\a^{3g-4} =\g\a^{3g-4} +s \a^{3g-3}$. Considering the component in
  $H^{6g-6}(\M)$ and using lemma~\ref{lem:8}, we have
  $$
   -2 \sum \Q_L^N(\a, \stackrel{(3g-4)}{\ldots}, \a, \q_i, \q_{i+g})=
   \Q_L^N(\a, \stackrel{(3g-4)}{\ldots}, \a, \g) + s <\a^{3g-3},[\M]>.
  $$
  Now, in $N$, $\g$ is the cup product $-2 \sum \p_i\p_{i+g} h^2$. It is easy
to 
  check that this coincides with the quantum product $-2 \sum \p_i\p_{i+g}
h^2$.
  For $g >3$ it is evident because of the degree. For $g=3$ there might be a
quantum
  correction in $H^0(N)$, but this is $-2 \sum \Q_L^N(\p_i,\p_{i+g},
h,h,\point)=0$
  (since lines are contained in the fibres).
  Now lemma~\ref{lem:10} and its proof imply that $-2 \sum \Q_L^N(\a, 
   \ldots, \a, \q_i, \q_{i+g})=\Q_L^N(\a, \ldots, \a, \g)$, so $s=0$.
\end{pf}

We are pursuing to prove an isomorphism between $QH^*(\M)$ and $HF^*(\Y)$, the

instanton Floer homology of the three manifold $\Y$.
First recall the main result contained in~\cite{Floer}.

\begin{thm}[\cite{Floer}]
\label{thm:15}
  Define $R^1_0=1$, $R^2_0=0$, $R^3_0=0$ and then recursively, for all 
  $r \geq 1$,
  $$
   \left\{ \begin{array}{l} R_{r+1}^1 = \a R_r^1 + r^2 R_r^2
   \\ R_{r+1}^2 = (\b+(-1)^{r+1}8) R_r^1 + {2r \over r+1}  R_r^3
   \\ R_{r+1}^3 =  \g  R_r^1
    \end{array} \right.
  $$
  Put $I'_r=(R^1_r,R^2_r,R^3_r) \subset \QQ[\a,\b,\g]$, $r \geq 0$. Then
  the $\Spz$-decomposition of $HF^*(\Y)$ is
$$
   HF^*(\Y)= \bigoplus_{k=0}^{g-1} \L_0^k H^3 \ox \QQ [\a, \b, \g]/I'_{g-k}.
$$
\end{thm}

The elements $R_r^1$, $R_r^2$ and $R_r^3$
are deformations graded mod $4$ of $q_r^1$, $q_r^2$ and $q_r^3$, respectively.
This means that we can write
\begin{equation}
\label{eqn:qu11}
   R_r^i= \sum_{j \geq 0} R_{r,j}^i,
\end{equation}
where $\deg(R_{r,j}^i)=\deg(q_r^i)-4j$, $j \geq 0$, and $R_{r,0}^i=q_r^i$.
In the case of $QH^*(\M)$ we shall have

\begin{prop}
\label{prop:16}
   The $\Spz$-decomposition of $QH^*(\M)$ is
   $$
   QH^*(\M)= \bigoplus_{k=0}^{g-1} \L_0^k H^3 \ox \QQ [\a, \b, \g]/J_{g-k},
   $$
   where $J_r$ is generated by three elements $Q_r^1$, $Q_r^2$ and $Q_r^3$,
which
   are deformations graded mod $4$ of $q_r^1$, $q_r^2$ and $q_r^3$,
respectively.  
\end{prop}

\begin{pf}
The action of $\Spz$ on $\M$ being symplectic (see~\cite[section
3.1]{Siebert}), we
have an epimorphism of rings (with $\Spz$-actions) like in~\eqref{eqn:qu3}
$$
  \AA(\S) \surj QH^*(\M).
$$
This induces an epimorphism on the invariant parts
$$
  \QQ[\a,\b,\g ]  \surj QH^*_I(\M),
$$
where $\g=-2 \sum_{i=0}^g \q_i\q_{i+g}$ (see lemma~\ref{lem:14}). Therefore
we have maps
\begin{equation}
  \L_0^k (\seq{\q}{1}{2g}) \ox \QQ[\a,\b,\g ]  \ar QH^*(\M).
\label{eqn:qu12}
\end{equation}
Let $V_k$ be the image of the map~\eqref{eqn:qu12}. As $\L_0^k H^3$, $0 \leq k
\leq g-1$,
are inequivalent irreducible $\Spz$-modules, the subspaces $V_k$ are pairwise
orthogonal.
On the other hand, the existence of the basis~\eqref{eqn:antes2} of $H^*(\M)$ 
and the results in~\cite{ST2} imply
that $\{x^{(k)}_i\a^a\b^b\g^c/ k=0,1,\ldots, g-1, \> a+b+c < g-k, \> i \in B_k
\}$ 
(where quantum
products are now understood) is a basis of $QH^*(\M)$. So the subspaces $V_k$
generate
$QH^*(\M)$, i.e.
\begin{equation}
  QH^*(\M)= \bigoplus_{k=0}^{g-1} V_k.
\label{eqn:qu13}
\end{equation}
Actually this decomposition coincides with the decomposition~\eqref{eqn:antes}.
This
is proved by giving a definition of $V_k$ independent of the ring structure
(cup 
product or quantum product). For instance, say that 
$V_k$ is the space generated by elements which are
orthogonal to $V_0$, $\ldots$, $V_{k-1}$ and such that the $\Spz$-module
generated by
them have dimension equal to $\dim \L_0^k H^3$.

Our second purpose is to describe the kernel 
of~\eqref{eqn:qu12}, i.e. the relations satisfied by the 
elements of $\L_0^k H^3$, $\a$, $\b$ and $\g$. The results in~\cite{ST2} imply
that 
we only need to write the relations of $V_k \subset
H^*(\M)$ in terms of the quantum product.
Fix $k$, and recall $x^{(k)}_0=\q_1\q_2\cdots \q_k \in \L_0^k H^3$
from~\eqref{eqn:antes3}. By
section~\ref{sec:ordinary} the relations in $V_k$ are given by $x^{(k)}_0
q_{g-k}^i$, $i=1,2,3$, and its $\Spz$-transforms.
We rewrite these relations in terms of the quantum product, using the 
basis~\eqref{eqn:antes2}, as
\begin{equation}
  x^{(k)}_0 q_{g-k}^i= \sum_{a+b+c<g-k} x_{abc} \a^a\b^b\g^c \in QH^*(\M)
\label{eqn:qu14}
\end{equation}
where $x_{abc} \in \L^k_0 H^3$, and the monomials in the right hand side have
degree
strictly less that the degree of the left hand side. Now we want to prove
that $x_{abc}$ are all multiples of $x^{(k)}_0$. Suppose not. Then it is easy
to see
that there exists $\p \in \Spz$ satisfying $\p( x^{(k)}_0)= x^{(k)}_0$ and
$\p(x_{abc}) \neq x_{abc}$. Consider~\eqref{eqn:qu14} minus its transform
under
$\p$. This is a relation between the elements of the basis of $V_k$, which is
impossible.

Therefore~\eqref{eqn:qu14} can be rewritten as
$$
  x^{(k)}_0 (q_{g-k}^i+  Q_{g-k,1}^i + Q_{g-k,2}^i +\cdots )=0,
$$
where $\deg Q_{g-k,j}^i= \deg q_{g-k}^i- 4j$, $j\geq 1$. This finishes the
proof.
\end{pf}

\begin{lem}
\label{lem:17}
  $\g J_k \subset J_{k+1} \subset J_k$, for $k=0,1,\ldots, g-1$.
\end{lem}

\begin{pf}
  Let $f\in J_k \subset \QQ[\a,\b,\g]$. By definition
(proposition~\ref{prop:16})
  this means that the quantum product $\q_1 \cdots \q_{g-k} f=0$. Using the
  action of $\Spz$ we have  
  $\q_1 \cdots \q_{g-k-1} \q_i f=0$, for $g-k \leq i \leq g$. 
  Thus $\q_1 \cdots
  \q_{g-(k+1)} \g f=0$, i.e. $\g f \in J_{k+1}$. For the second inclusion,
  let $f\in J_{k+1}$. Then $\q_1 \cdots \q_{g-(k+1)} f=0$ and hence 
  $\q_1 \cdots \q_{g-k} f=0$, i.e. $f \in J_k$.
\end{pf}

\begin{prop}
\label{prop:18}
  There are numbers $c_r, d_r \in \QQ$, $1 \leq r \leq g-1$, 
  such that for $0 \leq r \leq g-1$ it is
  $$
   \left\{ \begin{array}{l} Q_{r+1}^1 = \a Q_r^1 + r^2 Q_r^2
   \\ Q_{r+1}^2 = (\b+c_{r+1}) Q_r^1 + {2r \over r+1}  Q_r^3
   \\ Q_{r+1}^3 = \g  Q_r^1 + d_{r+1} Q_r^2
    \end{array} \right.
  $$
\end{prop}

\begin{pf}
  Completely analogous to the proof of~\cite[theorem 10]{Floer}.
\end{pf}

\begin{prop}
\label{prop:19}
  For all $1 \leq r \leq g$, 
  $c_r= (-1)^{r+g+1}\, 8$ and $d_r=0$.
\end{prop}

\begin{pf}
  We write $R_{g-k}^i= \sum_{j \geq 0} R_{g-k,j}^i$ and $Q_{g-k}^i= \sum_{j
\geq 0} 
  Q_{g-k,j}^i$, as in~\eqref{eqn:qu11}, for $0 \leq k \leq g-1$.
  Then $R_{g-k,0}^i=Q_{g-k,0}^i=q_{g-k}^i$. The coefficients $c_r$ and $d_r$
are
  determined by the first correction term $Q_{r,1}^i$ of $Q_r^i$. By the
definition
  of $R^i_r$ in theorem~\ref{thm:15}, we only need to check that 
  $R_{g-k,1}^i=(-1)^g Q_{g-k,1}^i$, for $i=1,2,3$, $0 \leq k \leq g-1$.

  Fix $i$ and $k$. Recall $x_0^{(k)}=\q_1 \cdots \q_k \in \L_0^k H^3$.
  By theorem~\ref{thm:15},
  $x_0^{(k)} R^i_{g-k} =0 \in HF^*(\Y)$. Pick an arbitrary $f=
  \a^a\b^b\q_{i_1}\cdots\q_{i_r} \in \AA(\S)$ of degree $6g-2-\deg(
  x_0^{(k)} q^i_{g-k})$. In $HF^*(\Y)$ the pairing 
  $<x_0^{(k)} R^i_{g-k},f>=0$, i.e.
  $D^{(w,\S)}_{S,H} (\bar x_0^{(k)} \bar R^i_{g-k} \bar f)=0$, where
  $\bar R^i_{g-k}=R^i_{g-k}(2\S,-4 x,-2\sum \g^{\#}_i\g^{\#}_{g+i})$,
  and analogously for $\bar f$ and $\bar x_0^{(k)}$
  (for the notation $D^{(w,\S)}$ see~\cite{Floer}). This means that
  $$
   D^{\PP^1}_{S,H} (\bar x_0^{(k)}  \bar R^i_{g-k,1} \bar f) 
  + D^{c_1}_{S,H} (\bar x_0^{(k)}  \bar R^i_{g-k,0} \bar f)=0.
  $$
  From theorem~\ref{thm:11} and remark~\ref{rem:13} we have that 
  the component in $H^{6g-6}(\M)$ of the quantum product
  $- x_0^{(k)}R^i_{g-k,1} f + (-1)^{g-1}
   x_0^{(k)}R^i_{g-k,0} f \in QH^*(\M)$ 
  vanishes, i.e.
  \begin{equation}
     -< x_0^{(k)}R^i_{g-k,1}, f> + (-1)^{g-1}<
    x_0^{(k)}R^i_{g-k,0} ,f>=0 \label{eqn:main1}
  \end{equation}
  in $QH^*(\M)$.

  On the other hand, proposition~\ref{prop:16} says that 
  $x_0^{(k)} Q^i_{g-k} =0 \in QH^*(\M)$. Multiplying by $f$,
  $x_0^{(k)}Q_{g-k}^i f =0$, so the component in 
  $H^{6g-6}(\M)$ of the quantum product $x_0^{(k)}Q^i_{g-k,1} f +
   x_0^{(k)}Q^i_{g-k,0} f \in QH^*(\M)$  is zero. Thus 
  \begin{equation}
          < x_0^{(k)}Q^i_{g-k,1}, f> + <
    x_0^{(k)}Q^i_{g-k,0} ,f>=0 \label{eqn:main2}
  \end{equation}
  in $QH^*(\M)$.

  Equations~\eqref{eqn:main1} and~\eqref{eqn:main2} imply together that
  \begin{equation}
  < x_0^{(k)}Q^i_{g-k,1}, f>= (-1)^g < x_0^{(k)}R^i_{g-k,1}, f>,
  \label{eqn:main3}
  \end{equation}
  for any $f \in \AA(\S)$ of degree $6g-2-\deg(
   x_0^{(k)}q^i_{g-k})= 6g-6-\deg( x_0^{(k)} Q^i_{g-k,1})$.
  As we are considering the pairing on classes of complementary degree, 
  equation~\eqref{eqn:main3} holds in  $H^*(\M)$ as well. 
  By section~\ref{sec:ordinary},
  $Q^i_{g-k,1} \equiv (-1)^g R^i_{g-k,1} \pmod{I_{g-k}}$.
  Considering the degrees, it must be $Q^1_{g-k,1} =(-1)^g R^1_{g-k,1}$
  and $Q^2_{g-k,1} =(-1)^g R^2_{g-k,1}$. For $i=3$, the difference
  $Q^3_{g-k,1} -(-1)^g R^3_{g-k,1}$ is a multiple of $q^1_{g-k}$. 
  The vanishing of the coefficient of $\a^{g-k}$
  for both $Q_{g-k}^3$ and $R_{g-k}^3$ (see theorem~\ref{thm:15} 
  and equation~\eqref{eqn:qu14}) implies $Q^3_{g-k,1} =(-1)^g R^3_{g-k,1}$.
\end{pf}

  Putting all together we have proved the following
\begin{thm}
\label{thm:main}
  The quantum cohomology of $\M$, 
  for $\S$ a Riemann surface of genus $g \geq 3$, has
  a presentation
  $$
   QH^*(\M)= \bigoplus_{k=0}^{g-1} \L_0^k H^3 \ox \QQ[\a,\b,\g] /J_{g-k}.
  $$
  where $J_r=(Q^1_r, Q^2_r,Q^3_r)$ and $Q^i_r$ are defined recursively
  by setting 
  $Q^1_0=1$, $Q^2_0=0$, $Q^3_0=0$ and putting for all $r \geq 0$
  $$
   \left\{ \begin{array}{l} Q_{r+1}^1 = \a Q_r^1 + r^2 Q_r^2
   \\ Q_{r+1}^2 = (\b+(-1)^{r+g+1}8) Q_r^1 + {2r \over r+1}  Q_r^3
   \\ Q_{r+1}^3 = \g  Q_r^1
    \end{array} \right.
  $$
\end{thm}

\begin{cor}
\label{cor:21}
  Let $\S$ be a Riemann surface of genus $g \geq 3$. Then there is an
isomorphism
  $$ 
    QH^*(\M) \isom HF^*(\Y).
  $$
  For $g$ even, the isomorphism sends $(\a,\b,\g) \mapsto (\a,\b,\g)$.
  For $g$ odd, the isomorphism sends $(\a,\b,\g) \mapsto 
  (\sqrt{-1}\, \a,-\b,-\sqrt{-1}\, \g)$.
\end{cor}

\begin{pf}
  This is a consequence of the descriptions of $QH^*(\M)$ and $HF^*(\Y)$ in
  theorem~\ref{thm:main} and theorem~\ref{thm:15}, respectively.
\end{pf}

\begin{rem}
\label{rem:22}
  Alternatively, we can say that for any $g \geq 1$ there is an isomorphism 
  $QH^*(\M) \isom HF^*(\Y)$, taking  $(\a,\b,\g) \mapsto 
  (\sqrt{-1}^{\,g}\a,\sqrt{-1}^{\,2g}\b,\sqrt{-1}^{\,3g}\g)$.
\end{rem}

\section{The cases $g=1$ and $g=2$}
\label{sec:6}

Let us review the cases of genus $g=1$ and $g=2$ in the view of
theorem~\ref{thm:main}.
These cases are somehow atypical, as the generators precise the introduction of

quantum corrections, a fact already noted in~\cite{Vafa}.

\begin{ex}
\label{ex:23}
  Let $\S$ be a Riemann surface of genus $g=1$. Then $\M$ is a point and we
can
  write
  $$
   QH^*(\M)= \QQ[\a,\hat \b,\g] /(\a,\hat \b+8,\g),
  $$
  where we have defined $\hat \b=\b-8$. This agrees with
theorem~\ref{thm:main}
  but with corrected generators. Again
  $QH^*(\M) \isom HF^*(\Y)$, where $(\a,\hat \b,\g) \mapsto 
  (\sqrt{-1}\, \a,-\b,-\sqrt{-1}\, \g)$.
\end{ex}

\begin{ex}
\label{ex:24}
  Let $\S$ be a Riemann surface of genus $g=2$.
  The quantum cohomology ring $QH^*(\M)$ has been computed by
  Donaldson~\cite{D1}, using an explicit description of $\M$ as the
  intersection of
  two quadrics in $\PP^5$. Let $h_2$, $h_4$ and $h_6$ be the integral 
  generators
  of $QH^2(\M)$, $QH^4(\M)$ and $QH^6(\M)$, respectively.
  Then, with our notations, $\a=h_2$, $\b = -4 h_4$ and
  $\g = 4 h_6$ (see~\cite{Vafa}). Define $\hat \g =-2 \sum \q_i\q_{i+g} \in
QH^*(\M)$.
  The computations in~\cite{D1} yield
  $\hat \g =\g -4 \a$ (compare with lemma~\ref{lem:14}). Put
  $\hat \b=\b + 4$. It is now
  easy to check that the relations found in~\cite{D1} can be 
  translated to 
  $$
   QH^*(\M)= \left( H^3 \ox \QQ[\a,\hat \b,\hat \g] /(\a,\hat \b-8,\hat \g)
\right)
   \oplus \QQ[\a,\hat \b,\hat \g] /(Q_2^1,Q_2^2,Q_2^3),
  $$
  where $Q_2^1= \a^2+\hat \b-8$, $Q_2^2= 
  (\hat \b+8) \a +\hat \g$ and $Q_2^3= \a \hat \g$ (defined exactly as in 
  theorem~\ref{thm:main}, but with corrected generators).
  Now
  $QH^*(\M) \isom HF^*(\Y)$, where $(\a,\hat \b,\hat \g) \mapsto 
  (\sqrt{-1}\, \a,-\b,-\sqrt{-1}\, \g)$.
\end{ex}

  The artificially introduced definition of $\hat \b$ is due to the same
phenomenon
  which
  causes the failure of lemma~\ref{lem:10} for $g=2$, i.e. the quantum 
  product $h^2$ differs from $\b$ in~\eqref{eqn:qu9}
  (defined with the cup product) because of a quantum correction in 
  $QH^*(N)$ which appears when $g=2$.

\end{document}